\documentstyle[12pt]{article}
\newcommand{\beq}{\begin{eqnarray}}
\newcommand{\eeq}{\end{eqnarray}}
\newcommand{\nbeq}{\begin{eqnarray*}}
\newcommand{\neeq}{\end{eqnarray*}}

\newcommand{\be}[1]{\begin{equation} \label{#1}}
\newcommand{\ee}{\end{equation}}

\newcommand{\D}{\displaystyle}

\def\th{\theta}

\begin{document}
\tolerance=500

$ $\vspace{1.3in}
\begin{center}
{\Large \bf Empirical Bayes estimators for the reproduction
parameter of Borel-Tanner distribution}

\vspace{0.5cm}{\bf George P. Yanev\footnote{Address: Mathematics
Dept.,
4202 E. Fowler Ave, PHY114, University of South Florida, Tampa, FL 33620, USA. E-mail: gyanev@cas.usf.edu.} \\
University of South Florida, Tampa, Florida}
\end{center}

\noindent {\bf ABSTRACT}\ We construct empirical Bayes estimators
for the reproduction parameter of Borel-Tanner distribution
assuming LINEX loss and prove their asymptotic optimality. Some
properties of the estimators regret risk are illustrated through
simulations.


\vspace{0.5cm}
\begin{center}
{\bf 1. INTRODUCTION}
\end{center}

\vspace{0.3cm}
 The probability mass function of the Borel-Tanner distribution is
\begin{equation} \label{bt2}
p(x|\theta, r)
     = a_r(x) \theta^{x-r}e^{-\theta x}, \qquad (x=r, r+1,
    \ldots)
\end{equation} where $0<\theta<1$, $r$ is a positive integer
and $ a_r(x)=rx^{x-r-1}/(x-r)! $

Initially (\ref{bt2}) was derived as the probability distribution
of the number of customers served in a queuing system. It also
appears in random trees and branching processes. More
specifically, it is the distribution of the total progeny in a
Galton-Watson process assuming Poisson reproduction, see Aldous
\cite{a} for recent applications. Our interest in estimating
$\theta$ stems from its role as reproduction number of an epidemic
infection modeled by a branching process, see Farrington et al.
\cite{FKG03}. We study nonparametric (with respect to the prior)
empirical Bayes (NPEB) estimators for $\th$. The NPEB estimation
procedures rely on the assumption for existence of a prior
distribution $G$ which, however, is unknown.  Consider independent
copies \ $ (X_1, \th_1), \ldots , (X_{n+1}, \th_{n+1}) $ \ of $(X,
\th)$, where $\th$ has a distribution $G$, and conditional on
$\th$, $X$ has a Borel--Tanner distribution given by (\ref{bt2}).
The ``past'' data consist of independent observations $x_1, x_2,
\ldots , x_n$ obtained with independent realizations $\th_1,
\th_2, \ldots , \th_n$ of $\th$, where the $X_i$s are observable
and the $\th_i$s are not observable. Denote by $\th_n(x)$ an
empirical Bayes estimator for $\th$ based on the ``past" data and
the ``present" observation $x_{n+1}=x$.
As Maritz and Lwin \cite{ml} point out, an advantage of using NPEB
estimators is the minimum assumptions on the class of prior
distributions. It turns out that in the case of Borel-Tanner
distribution the Bayes rule assuming LINEX loss depends on the
prior through the marginals only. This remarkable fact allows us
to construct simple NPEB estimators estimating the Bayes rule
directly without estimating the prior itself.

Notice that NPEB estimators for $\theta$ under weighted
squared-error loss are studied in Yanev \cite{y01}. In the next
section we use the asymmetric LINEX loss function, instead. In
Section~3 we prove the estimators asymptotic optimality. The last
section contains numerical results concerning the estimators
performance measured by their regret risk.

\vspace{0.5cm}
\begin{center}
{\bf 2. EMPIRICAL BAYES ESTIMATION USING LINEX LOSS}
\end{center}

\vspace{0.3cm} In some applications (e.g. surveillance of
infectious diseases) the squared-error loss function seems
inappropriate in that it assigns the same loss to overestimates as
to equal underestimates. A well-known alternative (see Huang et
al. \cite{HuLi97} and the references therein) is the LINEX loss
function defined, for $\gamma_l \leq \gamma \leq \gamma_u$ and
$\gamma \neq 0$ \ by \be{E:EBLinex} L_{\gamma}(\hat{\th}, \th)=
e^{\gamma(\hat{\th}-\th)}-\gamma(\hat{\th}-\th)-1 \ , \ee where
$\hat{\th}$ is an estimator for $\th$. It is clear that the LINEX
loss function is convex, asymmetric and for $\gamma >0$ it
increases almost linearly for negative errors and almost
exponentially for positive errors. Thus, it penalizes an
overestimation more seriously than an underestimation. This is
reversed when $\gamma<0$. For small values of $| \gamma |$ the
LINEX loss is close to the squared--error loss. From now on we
assume that $\gamma$ is a positive integer; the case $\gamma <0$
can be treated similarly.

Based on a single observation, the maximum likelihood estimator
$\th_{MLE}(x)$ for $\th$ is (e.g., Kumar \& Consul \cite{kc})
\begin{equation} \label{MLE} \th_{MLE}(x)=\frac{\D x - r}{\D
x}.\end{equation} Denote by ${\cal I}_{\D A}$ the indicator of the
event $A$.

{\bf Theorem 1}\ {\it Assume LINEX loss with $\gamma>0$, integer.
A NPEB estimator for $\theta$ in (\ref{bt2}) is
\begin{equation}\label{gammaestimate}
  \th_n(x)=
  \gamma^{-1}\ln \tau_n(x){\cal I}_{\D \left\{ \tau_n(x)\in (1,e^\gamma)\right\}} +
  (x-r)/x{\cal I}_{\D \left\{\tau_n(x)\notin (1,e^\gamma)\right\}},
\end{equation}
where
\[ \tau_n(x)= \frac{r+\gamma}{r}\left(\frac{\D x+\gamma}{\D x}\right)^{x-r-1}
       \! \!  \frac{\D m_n(x|r)}{\D m_n(x+\gamma|r+\gamma)}\
\]
and $m_n(z|y)$ is an estimate for the marginal distribution \
$m_G(z|y)=\int_0^1 p(z|\th, y)dG(\th)$. }

 \noindent {\bf Proof}\ \ The Bayesian estimator $\theta_{G}(x)$ under
LINEX loss is  (e.g. Huang et al. \cite{HuLi97}) \be{E:BELinex}
\th_G(x)=-\gamma^{-1}\ln E_{G|x} e^{-\gamma\th}, \ee provided that
$E_{G|x}e^{-\gamma\theta}<\infty$, where $E_{G|x}(\cdot)$ is the
expectation w.r.t. the posterior. Since
 \nbeq E_{G|x} e^{-\gamma \th}
        & = &
        \frac{\D 1}{\D m_G(x|r)}\int_0^{1} e^{-\gamma \th}a_r(x)\th^{x-r}e^{-x\th }dG(\th) \\
        & = &
        \frac{\D a_r(x)}
        {\D a_{r+\gamma}(x+\gamma)}\frac{\D m_G(x+\gamma|r+\gamma)}
        {\D m_G(x|r)},
\neeq
 we can write the Bayesian estimator $\th_G(x)$ from (\ref{E:BELinex}) as
 \beq
 \th_G(x) & = & \gamma^{-1}\ln
\left\{
        \frac{r+\gamma}{r}\left(\frac{\D x+\gamma}{\D x}\right)^{x-r-1}
        \frac{\D m_G(x|r)}{\D m_G(x+\gamma|r+\gamma)}
\right\} \nonumber \\
        & = & \gamma^{-1}\ln\tau_G(x), \quad \mbox{say.} \nonumber
\eeq
 Note that, $\th_G(x)$ depends upon the prior through the marginal
distribution only. Therefore, estimating the marginals, we can
construct a NPEB estimator $\theta_n(x)$ for $\th$ as given in
(\ref{gammaestimate}). \hfill{$\Box$}

One possible form of the estimators $m_n(z|y)$ in Theorem~1 can be
obtained as follows. In addition to the current $X_{n+1}(r)=x$,
let us have observed $n$
independent pairs 
\begin{equation}\label{scheme}
(X_1(r), X_1(\gamma)), (X_2(r), X_2(\gamma)), \ldots , (X_{n}(r),
X_{n}(\gamma)),
\end{equation}
where $X_i(r)$ and $X_i(\gamma)$ are independent and Borel-Tanner
distributed
 with $p(x|\th_i,r)$ and $p(x|\th_i,\gamma)$,
respectively. It is known (e.g., Kumar \& Consul \cite{kc}) that
$X_i(r)+X_i(\gamma)$ has pmf $p(x|\th_i, r+\gamma)$. Let
$f_{n}(y|r+\gamma)$ be the number of pairs,
 such that $X_i(r)+X_i(\gamma)=y$, $(i=1,\ldots, n)$. Consistent
estimators for the marginals $m_G(x+\gamma|r+\gamma)$ and
$m_G(x|r)$ are the relative frequencies \be{relfreq}
m_n(x+\gamma|r+\gamma)=\frac{f_{n}(x+\gamma|r+\gamma)}{n+1} \quad
\mbox {and} \quad
m_n(x|r)=\frac{1+f_n(x|r)}{n+1}.
\ee

Let us notice here that a NPEB estimator $\tilde{\theta}_n(x)$ for
$\theta$ under the squared-error loss $L(\hat{\theta},
\theta)=(\hat{\theta}-\theta)^2$ is constructed in Yanev
\cite{y01} as follows
\[
 \tilde{\th}_n(x)=
  \kappa_n(x){\cal I}_{\D \left\{ \kappa_n(x)\in (0,1)\right\}} +
  (x-r)/x{\cal I}_{\D \left\{\kappa_n(x)\notin (0,1)\right\}},
\]
where
\[ \kappa_n(x)= \frac{a_r(x)}{m_n(x)}\sum_{j=0}^\infty
\frac{(j+1)^{j-1}}{j!} \frac{m_n(x+j+1)}{a_r(x+j+1)},
\]
where, as before, $a_r(y)=ry^{y-r-1}/(y-r)!$
 \vspace{0.5cm}
\begin{center}
{\bf 3. ASYMPTOTIC OPTIMALITY }
\end{center}

\vspace{0.3cm} The Bayes risk of an estimator $\hat{\theta}$ can
be written as
\[R(G, \hat{\theta})=
\int_X\int_{\Theta}
            L(\hat{\theta}, \th)p(x|\theta, r) dG(\theta)dx \\
            =
\sum_{x=r}^\infty \int_{\Theta}
            L(\hat{\theta}, \th)p(\th|x) dG(\theta)m_G(x|r), \\
 \]

 \vspace{-0.4cm} \noindent where $p(\th |x)$ is the posterior
distribution. If $R(G, \th_n|\underline{X}_n)$ is the conditional
Bayes risk of the estimator $\th_n(x)$ given
$\underline{X}_n=(X_1, \ldots , X_n) $, then $ R(G, \th_n) =
E_n\{R(G, \th_n |\underline{X}_n)\} $ is the (unconditional) Bayes
risk of $\th_n$, where the expectation $E_n(\cdot)$ is taken with
respect to $\underline{X}_n$. The estimator $\th_n(x)$ is
asymptotically optimal for given $G$ if $ \lim_{n\to \infty}R(G,
\th_n)=R(G, \th_G). $ We shall prove the asymptotic optimality of
$\th_n(x)$.

First, let us find the minimum Bayes risk $R(G, \th_G)$ attained
by  the Bayesian estimator $\th_G(x)$. Since (\ref{E:BELinex})
implies $\exp\left(\gamma\th_G(x)\right)
\int_0^1\exp(-\gamma\th)p(\th |x)d\th=1$, we have
 \nbeq R(G,\th_G) & =
& \sum_{x=r}^{\infty}\left\{\int_0^1 \left\{ e^{\gamma
(\th_G(x)-\theta)}-\gamma(\th_G(x)-\theta)-1\right\} p(\th|x)d\th
\right\} m_G(x|r)
  \\
& = & \sum_{x=r}^{\infty}\left\{ e^{\gamma\th_G(x)}
    \int_0^1 e^{-\gamma\th}p(\th |x)d\th
    -\gamma\th_G(x) + \int_0^1\gamma\th p(\th|x)d\th-1\right\}
        m_G(x|r)
       \\
& = & \sum_{x=r}^{\infty}\left\{
        \int_0^1\gamma\th p(\th|x)d\th-\gamma \th_G(x)\right\}
        m_G(x|r).
\neeq Next, using $\exp\left(\gamma\th_G(x)\right)
\int_0^1\exp(-\gamma\th)p(\th |x)d\th=1$ again, we obtain
 \nbeq R(G, \th_n)
& = &  \sum_{x=r}^{\infty}E_n\left\{ e^{\gamma\th_n(x)}
    \int_0^1 e^{-\gamma\th}p(\th |x)d\th
    -\gamma\th_n(x) + \int_0^1\gamma\th p(\th|x)d\th-1\right\}
        m_G(x|r)
        \\
 & = & \sum_{x=r}^{\infty}E_n\left\{e^{\gamma
(\th_n(x)-\th_G(x))}
        -\gamma\th_n(x) +\int_0^1\gamma\th p(\th |x)d\th-1\right\}
        m_G(x|r)  .
  \neeq
Therefore, \be{risk} R(G,\th_n)-R(G, \th_G)
=\sum_{x=r}^{\infty}E_n\left\{e^{\gamma (\th_n(x)-\th_G(x))}
        -\gamma (\th_n(x)-\th_G(x)) -1\right\}
        m_G(x|r)
\ee

 Let us truncate the Borel--Tanner distribution (\ref{bt2}) starting with $r=k$ as
follows
\be{E:BTtrunc} p^\ast(x|\th, k)=
\left\{%
\begin{array}{ll}
    p(x|\th, k), & \mbox{if} \ k\leq x\leq k+N-1; \\
    \sum_{x=k+N}^\infty p(x|\th, k), & \mbox{if} \ x=k+N. \\
\end{array}%
\right.     \ee
 where $N$ is a positive integer. Denote the
truncated marginal by $m^\ast_G(x)=\int_0^1p^\ast(x|\th,
y)dG(\th)$. Similar to the non-truncated case, if $r\leq x\leq
r+N-1$ then
\[
E_{G|x} (e^{-\gamma \th})=
        \frac{\D a_r(x)}
        {\D a_{r+\gamma}(x+\gamma)}\frac{\D m^\ast_G(x+\gamma|r+\gamma)}
        {\D m^\ast_G(x|r)}=\frac{\D a_r(x)}
        {\D a_{r+\gamma}(x+\gamma)}\frac{\D m_G(x+\gamma|r+\gamma)}
        {\D m_G(x|r)}= \frac{\D 1}{\D \tau_G(x)}.
\]
If $x=r+N$ then
 \nbeq E_{G|x} (e^{-\gamma \th})
        & = &
        \frac{\D 1}{\D m^\ast_G(r+N|r)}
        \int_0^1e^{-\gamma \th}
        \sum_{k=r+N}^\infty a_r(k)\th^{k-r}e^{-\th k}dG(\th) \\
        & = &
        \frac{\D 1}{\D m^\ast_G(r+N|r)}
        \sum_{k=r+N}^\infty
         \frac{\D a_r(k)}{\D a_{r+\gamma}(k+\gamma)}
                 m^\ast_G(k+\gamma |r+\gamma)\\
    & = &
     \left( \sum_{k=r+N}^\infty
                 m_G(k|r)/\tau_G(k)\right)/\sum_{k=r+N}^\infty m_G(k|r).
\neeq  Let $\tau_G^\ast(x) =\tau_G(x)$ if $r\leq x\leq r+N-1$; \ \
$= \sum_{k=r+N}^\infty m_G(k|r)/\sum_{k=r+N}^\infty
                \left( m_G(k|r)/\tau_G(k)\right)$
                 if $x=r+N$.
The Bayesian estimator in the truncated case is given by $
 \th^\ast_G(x) =
\gamma^{-1}\ln\tau^\ast_G(x)$. Let us estimate $m_G(x|y)$ by
$m_n(x|y)$ as in Theorem 1 and set $\tau_n^\ast(x) =\tau_n(x)$ if
$r\leq x\leq r+N-1$; \ \ $= \sum_{k=r+N}^\infty
m_n(k|r)/\sum_{k=r+N}^\infty
                \left( m_n(k|r)/\tau_n(k)\right)$
                 if $x=r+N$. We construct a NPEB estimator
in the truncated case as follows
\[
\theta^\ast_n(x)=
  \gamma^{-1}\ln \tau_n^\ast(x){\cal I}_{\D \left\{ \tau_n^\ast(x)\in (1,e^\gamma)\right\}} +
  (x-r)/x{\cal I}_{\D \left\{\tau_n^\ast(x)\notin (1,e^\gamma)\right\}}.
\]
Now, we are in a position to prove the asymptotic optimality of
$\th_n(x)$.

{\bf Theorem 2}\  Assume prior $G$ with finite first moment. If
$m_n(z|y)$ is a consistent estimator for $m_G(z|y)$, then the NPEB
estimator $\th_n(x)$ given by (\ref{gammaestimate}) is
asymptotically optimal, i.e.,
\[
\lim_{n\to\infty} R(G, \th_n)=R(G, \th_G).
\]

\noindent {\bf Proof}\ \  Since $\th_G(x)$ is the Bayesian
estimator, we have $R(G,\th_n)>R(G, \th_G)$ and thus
\begin{equation}\label{E:first}
 \hspace{-0.5cm} R(G,\th_n)-R(G, \th_G)
\leq  |R(G,\th_n)-R(G, \th^\ast_n)|+|R(G,\th^\ast_n)-R(G,
\th^\ast_G)| +|R(G,\th^\ast_G)-R(G, \th_G)|
\end{equation}
To prove the theorem it is sufficient to show that the right hand
side of (\ref{E:first}) has $\lim_{N\to\infty}\limsup_{n}$ equals
zero, when $N$ is from (\ref{E:BTtrunc}). The truncated analog of
(\ref{risk}) leads to
\[
\hspace{-0.5cm} |R(G,\th^\ast_n)-R(G, \th^\ast_G)|
=\sum_{x=r}^{r+N}E_n\left\{e^{\gamma
(\th_n^\ast(x)-\th_G^\ast(x))}
        -\gamma (\th_n^\ast(x)-\th_G^\ast(x)) -1\right\}
        m_G(x|r)
\]
Since $m_n(z|y)$ is a consistent estimator for $m_G(z|y)$, we have
$\lim_{n\to\infty}\theta_n^\ast(x)=\theta_G^\ast(x), \ F^\infty
\mbox{- a.s.}$, where $F^\infty$ is the product measure induced by
$X_1, X_2, \ldots, X_n,\ldots$. Notice that, both $\th^\ast_n$ and
$\th^\ast_G$ are bounded. Indeed, $\th^\ast_n$ is bounded by
definition and
 $0< \theta^\ast_G(x)= -(1/\gamma)\ln E_{G|x} (e^{-\gamma \th})
<(1/\gamma)\ln e^{\gamma}=1$. Therefore, by the Lebesgue dominated
convergence theorem we can pass to the limit inside the
expectation in the right hand side above
and
obtain \be{E:third} \lim_{n\to\infty}|R(G,\th^\ast_n)-R(G,
\th^\ast_G)|=0. \ee Also, since $p^\ast(\th |x,r)=p(\th |x,r)$,
$m^\ast_G(x)=m_G(x)$ for $r\leq x\leq r+N-1$, and
$m^\ast_G(r+N)=\sum_{x=r+N}^\infty m_G(x)$ it is not difficult to
obtain
\[
|R(G,\th^\ast_G)-R(G, \th_G)|  =  \sum_{x=r+N}^\infty \left\{
\int_0^1\gamma \th \left(p^\ast(\th |r+N)-p(\th |x)\right)d\th -
\gamma \left(\th^\ast_G (r+N)-\th_G(x)\right) \right\} m_G(x).
\]
 Since $|p^\ast(\th |r+N)-p(\th |x)|<1$,
$|\th^\ast_G (r+N)-\th_G(x)|<1$, and $E\th<\infty$ we have
\be{E:fourth} \lim_{N\to\infty}|R(G,\th^\ast_G)-R(G, \th_G)|=0.
\ee Similar to (\ref{E:fourth}) one can prove that $
\lim_{N\to\infty}|R(G,\th_n)-R(G, \th^\ast_n)|=0$.  This along
with (\ref{E:first})-(\ref{E:fourth}) completes the proof.
\hfill{$\Box$}

\vspace{2.4cm}
\begin{center}
{\bf 4. NUMERICAL EXAMPLES}
\end{center}

\vspace{0.3cm} Using the notation introduced before
(\ref{relfreq}) we set
\[ \tau_n^f(x)=
   \frac{r+\gamma}{r}\left(\frac{\D x+\gamma}{\D x}\right)^{x-r-1}
        \! \! \! \frac{\D f_n(x|r)}{\D f_n(x+\gamma|r+\gamma)}\ .
\]
Let $A=\left\{\tau_n^f(x)\in (1,e^\gamma)\cap
f_n(x+\gamma|r+\gamma)\neq 0\right\}$ and $A^c$ be its complement.
Making use of the relative frequency estimators (\ref{relfreq})
consider $\th_n^f(x)$ to be defined by
\[
  \theta_n^f(x)=
  \gamma^{-1}\ln \tau_n^f(x){\cal I}_{\D A} +
(x-r)/x{\cal I}_{\D A^c}.
\]
That is, if $A$ occurs, then we estimate $\th$ by $\gamma^{-1}\ln
\tau_n^f(x)$; whereas if $A^c$ occurs then we use the MLE
(\ref{MLE}) for $\th$ instead.

A popular measure of the performance of one estimator
$\hat{\th}(x)$ is its regret risk $S(\hat{\th})=R(G, \hat{\th})-
R(G, \th_G)>0$. For our simulation study we take $r=5$, Uniform
$(0.5,1)$ prior and LINEX loss with $\gamma=3$. Then the minimum
Bayes risk attained by the Bayesian estimator
 \nbeq \th_{U}(x)
& = &
    \frac{\D 1}{\D 3}\ln \frac{\D \int_{0.5}^1 \th^{x-5}e^{-x\th} d\th}
    {\D \int_{0.5}^1 \th^{x-5}e^{-(x+3)\th} d\th}\ ,
\neeq is $R(U_{(0.5, 1)}, \th_U)=0.0622$.

In the empirical Bayes scheme (\ref{scheme}), let us set $n=50$.
Selecting 50 random values for $\theta_i \sim U_{(0.5,1)}, \
i=1,2,\ldots, 50$, we generate two sets of 50 branching processes
starting  with $r=5$ and $\gamma=3$ ancestors, respectively and
both having $Poisson(\theta_i), i=1,2, \ldots 50$ offspring
distributions. Notice that the total progeny of each process is a
realization of a Borel-Tanner $(\theta_i, \cdot)$ random variable.
Repeating the above procedure 100 times, we obtain 100 samples of
50 pairs Borel-Tanner observations, $(X_i(5), X_i(3)), \
i=1,2,\ldots, 50$. Each sample gives us a NPEB estimate
$\theta_{50}^f(x)$ with regret risk $S_i(\theta_{50}^f), \ i=1, 2,
\ldots, 100$. We estimate the regret risk $S(\theta_{50}^f)$ with
the average $\bar{S}(\theta_{50}^f)=\sum_{i=1}^{100}
S_i(\theta_{50}^f)/100$.

The above scheme is repeated with $n=75$ and $n=100$.  As an
illustration,  we present in Table~\ref{tab:3} results for one
sample with $n=100$ . For this particular sample,
$S_i(\theta_{100}^f)=0.0980$, which is less than
$S(\theta_{MLE})=0.1327$.

\begin{table}[h]
\begin{center}
\vspace{0.5cm}
\begin{tabular}{|c|c|c|c|c|c|c|c|c|c|c|c|c|c|c|c|c|c|} \hline
  $x$
  & 5 & 6 & 7  & 8 & 9 & 10 & 11 & 12 & 13 & 14 & 15 & 16 & 17 &
  18 & 19 & 20
\\ \hline
 $\th_{100}^f(x)$
 & .46 &  .69  &  .92  &  .65  &  .58 &   .51   & .55  &  .53 & .62 & .96 & .61 & .69
& .16  &  .72  &  .79  &  .75
\\ \hline
  $\th_{U}(x)$ & {\it .63} &  {\it  .64} & {\it  .65} &  {\it .65} & {\it  .66} & {\it  .67} &
{\it .67}  & {\it .68} & {\it .69} & {\it .69} & {\it .70} & {\it
.71} & {\it .71} & {\it .72} & {\it .73} & {\it .73}
\\ \hline
 $\th_{MLE}(x)$
& 0  &  .16 &   .28 &   .38  &  .44   & .50    & .55 &   .58 & .62
& .64 &  .67 &   .69 &   .71 &   .72   & .74   & .75
\\ \hline
\end{tabular}
\caption{Estimates $\theta_n^f(x)$, $\th_{U}(x)$,  and
$\theta_{MLE}(x)$ for $\theta$ from a sample with
$n=100$.}\label{tab:3}
\end{center}
\end{table}

The numerical results for the regret risks are
given in Table~\ref{tab:1}. Several comments are in place. For
small $x$, (columns 2-4) and $n=75$ or 100, the improvement of
$\theta_n^f$ over $\theta_{MLE}$ is substantial. Overall, (columns
5-7), the regret risk of $\theta_n^f$ is not higher than that of
$\theta_{MLE}$.

\begin{table}[h]
\begin{center}
\vspace{0.3cm}
\begin{tabular}{|c||c|c|c||c|c|c|} \hline
& \multicolumn{3}{c||}{$0\leq x \leq 15$} &
\multicolumn{3}{c|}{$0\leq x \leq 200$} \\ \cline{2-7}
  $n$ & $\scriptstyle \bar{S}(\theta_n^f)$ & $\scriptstyle STD(\bar{S}(\theta_n^f))$  & $\scriptstyle S(\theta_{MLE})$
  & $\scriptstyle \bar{S}(\theta_n^f)$ & $\scriptstyle STD(\bar{S}(\theta_n^f))$  & $\scriptstyle S(\theta_{MLE})$ \\
  \hline \hline
  50 & 0.1211  &  0.0037 & 0.1292 & 0.1397 & 0.0037 & 0.1327  \\ \hline
  75 & 0.1076 & 0.0036 & 0.1292 & 0.1300 &  0.0037 & 0.1327   \\ \hline
  100 & 0.1038 &  0.0033  & 0.1292 & 0.1299 &  0.0036  & 0.1327 \\ \hline
\end{tabular}
\caption{Numerical results on regret risks of $\theta_n^f(x)$ and
$\theta_{MLE}(x)$.}\label{tab:1}
\end{center}
\end{table}

Finally, note that the Borel-Tanner distribution (\ref{bt2}) has
monotone likelihood ratio in $x$, i.e., $p(x|\th ', r)/p(x|\theta,
r)$ is an increasing function of $x$ whenever $0<\th<\th '<1$.
This suggests that the NPEB $\th_n(x)$ can be improved on by the
monotonizing procedure of Van Houwelingen and Stijnen \cite{hs93}.

\begin{center}
{\bf ACKNOWLEDGEMENTS}
\end{center}
I thank R. Gueorguieva  for helping me with the simulations done
with an Ox version 3.30. I also thank the referee for the valuable
comments. This research is partially supported by NFSI-Bulgaria,
Grant No. MM-1101/2001.



\begin{thebibliography}{99}

\bibitem{a}
Aldous, D.J. Deterministic and stochastic models for coalescence
(aggregation and coagulation): a review of the mean-field theory
for probabilists. Bernoulli {\bf 1999}, {\it 5}, 3-48.



\bibitem{FKG03}
Farrington, C.P.; Kanaan, C.P.; Gay, N.J. Branching process models
for surveillance of infectious diseases controlled by mass
vaccination. Biostatistics {\bf 2003}, {\it 4}(2), 279-295.



\bibitem{HuLi97}
Huang, S.Y.; Liang, T.C. Empirical Bayes estimation of the
truncation parameter with Linex loss. Statist. Sinica {\bf 1997},
{\it 7}, 755--769.



\bibitem{kc}
A. Kumar, A.; Consul, P.C. Minimum variance unbiased estimation
for modified power series distribution. Comm. Statist. A - Theory
Methods {\bf 1980}, {\it 9}, 1261-1275.



\bibitem{ml}
Maritz, J.S.; Lwin, T. {\it Empirical Bayes Methods}, 2nd Ed.;
Chapman and Hall: London, 1989.



\bibitem{hs93} Van Houwelingen, J.C.; Stijnen, T. Monotone
empirical Bayes estimators based on more informative samples. J.
Amer. Statist. Assoc. {\bf 1993}, {\it 88}, 1438-1443.

\bibitem{y01}
Yanev, G.P. {\it Statistical modeling of epidemic disease
propagation via branching processes and Bayesian inference},
Dissertation, 2001, University of South Florida.




\end{thebibliography}
\end{document}